\documentclass[12pt, a4paper]{amsart}
\usepackage{amsfonts,amsmath, amsthm, amssymb, amscd}
\usepackage{latexsym}
\usepackage{hyperref}

\def\goth#1{\mathfrak{#1}}

\newcounter{Theorem}
\newtheorem{theorem}[Theorem]{\sc Theorem}
\newtheorem{prop}[Theorem]{\sc{Proposition}}
\newtheorem{lemma}[Theorem]{\sc{Lemma}}
\newtheorem{cor}[Theorem]{\sc{Corollary}}

\newtheorem{dfn}[Theorem]{\sc{Definition}}
\newtheorem{example}[Theorem]{\sc{Example}}
\newtheorem{examples}[Theorem]{\sc{Examples}}

\newcommand {\PPP}{{{\sc{Proof of Theorems 1 and 2}$\ \ $}}}
\newcommand {\PPPP}{{{\sc{Proof of Theorem 3}$\ \ $}}}

\newcommand {\N}{\mathbb{N}} 
\newcommand {\R}{\mathbb{R}} 
\newcommand {\C}{\mathbb{C}} 

\newcommand {\A}{\mathcal{A}}
\newcommand {\B}{\mathcal{B}}

\newcommand {\D}{\mathcal{D}}
\newcommand {\E}{\mathcal{E}}
\newcommand {\bE}{\mathbf{E}}

\newcommand {\F}{\mathcal{F}}

\newcommand {\Pp}{\mathcal{P}}

\newcommand {\V}{\mathcal{V}}

\newcommand {\z}{\goth{z}}
\newcommand {\n}{\goth{n}}
\newcommand {\vv}{\goth{v}}

\DeclareMathOperator{\vol}{vol}

\def\qmbox#1{\quad\mbox{#1}\quad}

\begin{document}

\title[Spherical spectral synthesis on Damek-Ricci Spaces]{Spherical
  spectral synthesis and two-radius theorems on Damek-Ricci Spaces}
\author{Norbert Peyerimhoff * \and Evangelia Samiou **} 
\thanks{* Research supported by the University of Cyprus} 
\thanks{** Research supported by Alan Richards Fellowship at Grey
  College and Blaise Pascal Award} 
\date{29 February, 2008} 
\address{Department of Mathematical Sciences, University of Durham,
  Durham DH1 3LE, United Kingdom} 
\address{Department of Mathematics and Statistics, University of
  Cyprus, P.O. Box 20537, 1678 Nicosia, Cyprus}

\email{norbert.peyerimhoff@durham.ac.uk} 
\email{samiou@ucy.ac.cy}


\begin{abstract}
  We prove that spherical spectral analysis and synthesis
  hold in Damek-Ricci spaces and derive two-radius theorems.
\end{abstract}


\maketitle


\section{Introduction}

D. Pompeiu investigated in 1929 the following problem: Given a closed
set $K \subset \R^2$ of positive volume and a continuous function $f$
satisfying
\begin{equation} \label{pomprop}
\int_{\sigma(K)} f = 0
\end{equation}
for all rigid motions $\sigma$ of the plane. Does this imply that $f = 0$?

The question has a negative answer in the case of a closed disk. In
fact, for every radius $r > 0$ there are vectors $\xi \in \R^2$ such
that the functions
\[ f(x) = e^{i \langle x,\xi \rangle} \] have vanishing integrals over
all closed disks of radius $r > 0$. A longstanding question is whether
disks are essentially the only simply connected bounded
counterexamples (Pompeiu's problem). This is closely related to
Schiffer's conjecture (for more details see, e.g., \cite[Problem
80]{Y-82} and \cite{BST-73,Wi-81,K-93,CKS-00}).

However, if \eqref{pomprop} holds true for two different disks $K =
B_{r_1}$ and $K = B_{r_2}$ with radii $r_1, r_2 >0$ whose quotient
avoids a particular exceptional set, then $f \equiv 0$. This
exceptional set is given by the quotient of any two different positive
zeros of a particular Bessel function expression. For more details and
more general two-radius theorems we refer to \cite{Za-72}. Two-radius
problems have also been considered in more general geometries.  The
paper \cite{BZ-80} generalizes the above result to all rank
one symmetric spaces. In fact, for compact symmetric spaces of rank
one it is sufficient to check vanishing of integrals over all balls of
a single radius to conclude that $f \equiv 0$, as long as this radius
is not a zero of a particular Jacobi polynomial expression.

In this paper we prove two-radius results for Damek-Ricci
spaces. Damek-Ricci spaces are harmonic manifolds and comprise the
rank one symmetric spaces of noncompact type. Recently, Heber
\cite{Heb-05} proved that the non-flat simply connected {\em
  homogeneous} harmonic spaces are precisely the symmetric spaces of
rank one and the non-symmetric Damek-Ricci spaces.

To state our main results we first have to introduce some
notation. For a detailed discussion of the geometry and analysis of
Damek-Ricci spaces we refer the reader to the exposition
\cite{Rou-03}. A Damek-Ricci space is a semidirect product $X =
\R\ltimes N$ of a generalized Heisenberg group $N$ with $\R$. Let $\n$
be the Lie algebra of $N$ with $\z = [\n,\n]$ and $\vv = \z^\bot$ with
respect to the inner product of $\n$. We denote the dimensions of
$\vv$ and $\z$ by $p$ and $q$. Since the pair $(p,q)$ plays an
inportant role in our results we often write $X^{(p,q)}$ for the
Damek-Ricci space. In this context, we have $\R = X^{(0,0)}$. A
Damek-Ricci space $X$ is a solvable group carrying a non-positively
curved left invariant Riemannian metric.

In the sequel the {\em spherical functions} $\varphi_\lambda^{(p,q)}$
defined in terms of the hypergeometric function $F$ by
\begin{equation} \label{philampq} 
\varphi_\lambda^{(p,q)}(r) = F(\rho-i\lambda,\rho+i\lambda,n/2,
-\sinh^2(r/2)), 
\end{equation}
will be crucial. Here $n = p+q+1$ and $\rho = \frac{p}{4} +
\frac{q}{2}$. These are the radial eigenfunctions $\varphi_\lambda$ of
the Laplacian of $X = X^{(p,q)}$, see \eqref{radphi9}.

The space of smooth radial functions on $X$ endowed with the topology
of uniform convergence on compacta is denoted by $\E_0(X)$. A variety
$V$ is a proper closed subspace of $\E_0(X)$ which is invariant under
convolution with radial distributions of compact support.

\medskip

Now we can state the two-radius theorems for Damek-Ricci spaces. In
the sequel integration will always be with respect to the left
invariant metric on $X$. Let $B_r(x)$, $S_r(x)$ denote the geodesic
ball, geodesic sphere around $x\in X$ of radius $r$, respectively.
Spheres and balls of radius $r$ around the identity element $e \in X$
will be denoted by $S_r$ and $B_r$.

\begin{theorem} \label{theomain1}
  Let $X = X^{(p,q)}$ be a Damek-Ricci space and let $r_1, r_2 > 0$ be
  such that the equations
  \[ \varphi_\lambda^{(p,q+2)}(r_j)=0,\, j=1,2, \] 
  have no common solution $\lambda \in \C$.

  Suppose $f \in C(X)$ and
  \[ \int_{B_r(x)} f = 0 \] 
  for $r=r_1,r_2$ and all $x \in X$.  Then $f=0$.
\end{theorem}

Similiarly, for spherical averages, we have

\begin{theorem} \label{theomain2} 
  Let $X = X^{(p,q)}$ be a Damek-Ricci space and let $r_1, r_2 > 0$ be
  such that the equations
  \[ \varphi_\lambda^{(p,q)}(r_j)=0,\, j=1,2, \] 
  have no common solution $\lambda \in \C$.

  Suppose $f \in C(X)$ and
  \[ \int_{S_r(x)} f = 0 \] 
  for $r=r_1,r_2$ and all $x \in X$.  Then $f=0$.
\end{theorem}

Note that in Damek-Ricci spaces a function is harmonic if and only if
it satisfies the mean value property for all radii. In fact, it
suffices to have the mean value property for only two suitably chosen
radii in order to conclude harmonicity of a function:

\begin{theorem} \label{theomain3}
  Let $X = X^{(p,q)}$ be a Damek Ricci space and let $r_1, r_2 > 0$ be
  such that the equations
  \[ \varphi_\lambda^{(p,q)}(r_j)=1,\, j=1,2, \] 
  have no common solution $\lambda \in \C \backslash \{ \pm i \rho
  \}$.

  Then $f \in C^\infty(X)$ is harmonic if and only if
  \[ \frac{1}{{\rm vol}(S_r(x))} \int_{S_r(x)} f = f(x)\] 
  for $r=r_1,r_2$ and all $x \in X$.
\end{theorem}

Two-radius results are closely related to spectral analysis and
synthesis of the underlying space (see, e.g., \cite{BST-73}).
L. Schwartz \cite{Schw-47} proved that spectral synthesis holds on the
real line. Theorem \ref{blackydrintro} below carries over this result
to radial functions in Damek-Ricci spaces. For symmetric spaces of
rank one, spherical spectral synthesis was proved in \cite{BS-79}.
For further results on spectral synthesis in symmetric spaces see,
e.g., \cite{BG-86,Be-87,Wa-87}.

\begin{theorem}[Spherical Spectral Synthesis] \label{blackydrintro} 
  Let $X$ be a Damek Ricci space and $V$ a variety of
  radial functions. Then $V$ is the closure of the span of all functions
  $\varphi_{\lambda, k}=\frac{d^k}{d\lambda^k}\varphi_\lambda$ contained in
  $V$.
\end{theorem}

The article is organized as follows: In Section \ref{prelim2} we
introduce some basic properties and notions which are used throughout
this paper and discuss Schwartz's fundamental result. In Section
\ref{abel}, we introduce the Abel transform, prove some useful
properties and obtain a Paley-Wiener Theorem. These properties will be
used in Section \ref{sss} to derive spherical spectral synthesis for
Damek-Ricci spaces. Finally, Section \ref{tworad} is devoted to the
proofs of the above two-radius results.

\noindent {\bf Acknowledgements:} 
The first author is grateful for the hospitality of the University of
Cyprus. The second author thanks the Department of Mathematical
Sciences of Durham University and Grey College for the
hospitality. The second author is indebted to Mrs. Marjorie Roberts
and Dr. Roy Roberts.

\section{Preliminaries}\label{prelim2}

\subsection{The spaces $\D(X), \E(X), \D_0(X)$ and  $\E_0(X)$}

A Damek-Ricci space $X := X^{(p,q)}$ of dimension $n = p+q+1$ is a
semidirect product $\R \ltimes N$ where $N$ is a generalized
Heisenberg group of dimension $p+q$ with $q$-dimensional center. We
may thus write elements of $X$ as pairs $x = (t(x),n(x)) = t(x) \cdot
n(x)$. Note that $t: X \to \R$ is a group homomorphism.

By $\E(X)$ we denote the space of all smooth functions on $X$ with the
topology determined by the seminorms
\[ \Vert f \Vert_{D,K} = \sup_{x \in K} \vert D f(x) \vert, \] 
where $D$ is an arbitrary differential operator on $X$ and $K \subset
X$ is an arbitrary compact subset. $\E(X)$ is a Fr{\'e}chet space.

\begin{dfn}
  Let $r(x) = d(x,e)$ denote the distance of $x \in X$ from the
  identity $e$. The {\em averaging projector} $\pi:\E(X) \to \E(X)$ is
  defined by
  \[ \pi f(x) = \frac{1}{{\rm vol}\, S_{r(x)}} \int_{S_{r(x)}} f \
  . \]
\end{dfn}

Let $\E_0(X) = \pi \E(X)$ denote the space of all smooth radial
functions on $X$, equipped with the relative topology. The spaces
$\E(\R)$ and $\E_0(\R)$ are analogously defined. Note that $Q:
\E_0(\R) \to \E_0(X)$, $Q f(x) = f(r(x))$ is a topological
isomorphism.

Let $\D(X)$ denote the space of all smooth function on $X$ with
compact support and $\D_K(M)$, for $K \subset X$ compact, the subspace
of $\D(X)$ of functions with support in $K$. The topology of $\D(X)$
is the inductive limit topology of the spaces $\D_K(X)$, and $\D_K(X)$
has the induced topology of $\E(X)$ (see, e.g. \cite{Hel-84}). Again,
we have $\D_0(X) = \pi \D(X)$.

The convolution of $f \in\E(X)$ and $g \in\D(X)$ (or $f \in\D(X)$ and
$g \in\E(X)$) is defined as
\[ 
f*g(x) = \int_X f(y) g(y^{-1}x) \ dy= \langle f, (\check g)_{x^{-1}} \rangle, 
\]
where $\check g(x) := g(x^{-1})$ and $g_x(y) := g(xy)$.

We quote some useful properties of $\pi$ from \cite{DR-92}:
\begin{eqnarray}
\pi^2 &=& \pi, \\
\langle \pi f, g \rangle &=& \langle f, \pi g \rangle, \\
\pi (f * \pi g) &=& \pi f * \pi g. \label{piprop}
\end{eqnarray}
The convolution in $\E(X)$ is associative but not
commutative. However, property \eqref{piprop} yields commutativity of
the convolution for radial functions and $\D_0(X) *
\E_0(X)=\E_0(X) * \D_0(X) \subset \E_0(X)$.

Let $\Delta = {\rm div}\, {\rm grad}$ denote the Laplacian on $X$. For
every $\lambda \in \C$ there exists a unique radial function
$\varphi_\lambda \in \E_0(X)$ satisfying
\[ 
\Delta \varphi_\lambda = - (\lambda^2 + \rho^2) \varphi_\lambda, \quad 
\varphi_\lambda(e) = 1, 
\]
where $\rho = \frac{p}{4} + \frac{q}{2}$ (see \cite{Rou-03}). We have
\begin{equation}\label{radphi9} \varphi_\lambda(x) =
  \varphi_\lambda^{(p,q)}(r(x)) \ , \end{equation}
where $\varphi_\lambda^{(p,q)}$ was defined in \eqref{philampq}. 

\begin{remark}
The parameter $\rho= \frac{p}{4} + \frac{q}{2}$ of the Damek-Ricci 
space $X = X^{(p,q)}$ can be interpreted isoperimetrically, asymptotically and 
in terms of the spectrum:

\noindent
$\bullet$ We have
\[
  h(X) := \inf_{K \subset X\, {\rm compact}}
  \frac{{\rm area}(\partial K)}{{\rm vol}(K)} = 2 \rho.
\]
This follows from \cite[p. 66]{Rou-03} and the explicit {\em isoperimetric
Cheeger constant} calculation in \cite{PS-03}.   

\noindent
$\bullet$ Using ${\rm vol}(S_r) = \frac{2^n \pi^{n/2}}{\Gamma(n/2)} 
\sinh^{p+q}(r/2) \cosh^q(r/2)$, one easily verifies the identity
\[ 
\lim_{r \to \infty} \frac{\log {\rm vol}(S_r)}{r} = 2 \rho
\]
for the {\em logarithmic volume growth} of spheres. 

\noindent
$\bullet$ We have
\[ \sigma(\Delta) = (-\infty,-\rho^2] \] for the {\em spectrum of the
Laplacian} in the Hilbert space $L^2(X,\mu)$.  This follows from the
fact that the spherical Fourier transform $\F f(\lambda) = \langle f,
\varphi_\lambda \rangle$ on $\D_0(X)$ extends to a Hilbert space
isomorphism (see \cite[Thm 15]{Rou-03}) and that $\Delta$ transforms
under this isomorphism into the multiplication operator $g \mapsto
-(\lambda^2 + \rho^2) g$.

\end{remark}

\subsection{The spaces $\D'(X), \E'(X), \D_0'(X)$ and  $\E_0'(X)$}

We denote by $\E'(X)$ the dual of $\E(X)$, endowed with the strong
dual topology.  $\E'(X)$ is the space of distributions of compact
support on $X$. The spaces $\D'(X), \D_0'(X)$ and $\E_0'(X)$ are
defined analogously.

The convolution of two distributions $S \in \D'(X)$, $T \in \E'(X)$ (or
$T \in \D'(X)$, $S \in \E'(X)$) can be calculated as follows
\[ 
\langle S*T, f \rangle = \langle S, x \mapsto \langle T, y \mapsto
f(xy) \rangle \rangle.
\]
The space $\D(X)$ is contained in $\E'(X)$ via $f \mapsto T_f$,
$\langle T_f, g \rangle := \langle f, g \rangle$ and we have $T * T_f
= T_{T * f}$ with $(T * f)(x) = \langle T, (\check f)_{x^{-1}}
\rangle$. Using a Dirac sequence $\rho_\epsilon \in \D_0(X)$ we have,
for $T \in \E'(X)$, $T * \rho_\epsilon \to T$ for $\epsilon \to 0$,
which shows that $\D(X)$ is dense in $\E'(X)$.  Therefore, all above
properties for functions carry over to distributions.

The spherical Fourier transform of a distribution $T \in \E_0'(X)$ is
defined as
\[ \F T(\lambda) = \langle T, \varphi_\lambda \rangle. \] 
If $T \in \E_0'(\R)$, we have for the classical Fourier transform
\[ 
\hat T(\lambda) = \langle T, \phi_\lambda \rangle = \langle T,
\psi_\lambda \rangle,
\]
where $\phi_\lambda(t) = e^{i\lambda t}$ and $\psi_\lambda(t) =
\frac{1}{2} (e^{i\lambda t} + e^{-i\lambda t})$.

\subsection{Schwartz's result}

Mean periodic functions were first introduced and studied by Delsarte
in a paper of 1935. L. Schwartz \cite{Schw-47} proposed the following
intrinsic definition of mean periodic functions: a function $f \in
\E(\R)$ is called {\em mean periodic} if not every function in
$\E(\R)$ can be obtained as a limit of finite linear combinations of
translates $f_x(y) = f(x+y)$ of $f$. The vector space of functions
obtained as such limits is called the variety $V^f$ of $f$.
Equivalently, $V^f$ can be defined as the closure of all functions of
the type $T * f$ with $T \in \E'(\R)$. Thus, a non-zero function $f
\in \E(\R)$ is  mean periodic if $V^f$ is a proper subspace of
$\E(\R)$.

We denote by
\[ \phi_{\lambda,k}(t) := \frac{d^k}{d\lambda^k} e^{i\lambda t} =
i^k t^k e^{i\lambda t}. \] 
The spectrum ${\rm spec}\, f$ of $f$ is defined as follows:
\[ {\rm spec}\, f := \{ \phi_{\lambda,k} \in V^f \mid k \in \N_0,
\lambda \in \C \}. \]

L. Schwartz proved the following fundamental result:

\begin{theorem} \label{blackystrong}
Let $f \in \E(\R)$ be a non-zero mean periodic function. Then $f$ 
is the limit of finite linear combinations of functions in ${\rm spec}\, f$.
\end{theorem}

Schwartz's result actually means that {\em spectral synthesis} holds
in $\E(\R)$. {\em Spectral analysis} is the weaker statement that
${\rm spec}\, V^f$ is not empty for every non-zero mean periodic function $f$.

\medskip

We now adapt the above theorem to $\E_0(\R)$, the subspace of $\E(\R)$
of even functions on $\R$. The averaging projector $\pi: \E(\R) \to
\E_0(\R)$ is the canonical projection $(\pi f)(x):= \frac{1}{2}(f(x) +
f(-x))$. Let
\[ 
\psi_{\lambda,k}(t) := \frac{1}{2}(\phi_{\lambda,k}(t) + \phi_{\lambda,k}(-t)). 
\]

\begin{theorem} \label{blackyeven} 
  Let $f \in \E_0(\R)$ be a non-zero mean periodic function, i.e.,
  \[
  \{ 0 \} \neq V^f_0 := \overline{\{ T * f \mid T \in \E_0'(\R) \}}
  \neq \E_0(\R).
  \]
  Then $f$ is the limit of finite linear combinations of functions in
  \[ {\rm spec}_0\, f := \{ \psi_{\lambda,k} \in V^f_0\backslash \{ 0 \} \mid
  k \in \N_0, \lambda \in \C \}.
  \]
\end{theorem}

\begin{proof}
  As a consequence of \eqref{piprop} and continuity of $\pi$, we
  obtain for a mean periodic function $f \in \E_0(\R)$ that
  \begin{equation} \label{vv0} V^f = \overline { \{ T * f \mid T \in
      \E'(\R) \} } \subset \pi^{-1}(V^f_0) \neq \E(\R).
  \end{equation}
  Consequently, $f$ is also mean periodic in $\E(\R)$.  Now using
  Theorem \ref{blackystrong}, $f$ is the limit of functions $f_j$
  which are finite linear combinations of functions in ${\rm spec}\,
  f$. Then we also have
  \[ \pi(f_j) \to \pi(f) = f. \] 
  Since \eqref{vv0} implies that $\pi({\rm spec}\, f) \subset {\rm
    spec}_0\, f \cup \{ 0 \}$, $\pi(f_j)$ is a finite linear combination of
  functions in ${\rm spec}_0 \, f$, finishing the proof. \qed
\end{proof}


\section{The Abel transform on distributions}\label{abel}

\subsection{The Abel transform}

The Abel transform will be of great importance in our considerations.

\begin{dfn}\label{abtra}
  Let $j$ and $a$ be the maps
  $$j\colon\E_0(\R)\to\E(X)\qmbox{with} jf(x)=e^{\rho t(x)} f(t(x))$$
  and
  $$
  a\colon\E_0(\R)\to\E_0(X)\qmbox{with} a=\pi\circ j,
  $$
  i.e.,
  $$
  af(x)=\frac{1}{\vol S_{r(x)}}\int_{S_{r(x)}}e^{\rho t(y)}
  f(t(y))\ .
  $$
  The {\em Abel transform} $\A$ is then defined as the dual of $a$,
  i.e., as the map
  $$\A\colon\E_0'(X)\to\E_0'(\R)\qmbox{with}\langle\A T,f\rangle=\langle T, 
  af\rangle $$ 
  for distributions $T\in\E_0'(X)$ of compact support and smooth
  functions $f\in\E_0(\R)$.
\end{dfn}

\begin{remarks}

\noindent
$\bullet$ The restriction of $\A$ to $\D_0(X) \subset \E_0'(X)$ is
explicitly given by (see \cite{Rou-03}):
$$\A f(t) = e^{\rho t} \int_N f(t n).$$

\noindent
$\bullet$ We have
\begin{equation}\label{abelfilambda}
a\psi_\lambda=\varphi_\lambda \qmbox{and} a\psi_{\lambda,k}=\varphi_{\lambda,k}.
\end{equation}
For the first equation, see \cite[p. 80]{Rou-03}. The
second equation is obtained by differentiating this with respect to
$\lambda$. The spherical Fourier transform can be expressed in terms
of the Abel transforms by
\begin{equation}\label{fourabel}
\F T(\lambda)=\langle T,\varphi_\lambda\rangle=\langle T, a\psi_\lambda\rangle=
\langle\A T,\psi_\lambda\rangle=\widehat{\A T}(\lambda) \\
\end{equation}
for $T\in\E_0'(X)$.

\end{remarks}

\subsection{Properties of the Abel transform}

A key property of the Abel transform is that it preserves the
convolution.

\begin{prop}
For $T,S\in\E_0'(X)$ and $f\in\E_0(\R)$ we have
\begin{eqnarray}
\label{abelalg}\A (T*_XS) &=& \A T *_\R \A S, \\
\label{abelalg2} a(\A T*_\R f) &=& T *_X af.
\end{eqnarray}
\end{prop}

\begin{proof} 
  Recall that for radial distributions,
  $$
  \langle T*S, \phi\rangle=\langle T, x\mapsto \langle S, y\mapsto
  \phi(xy)\rangle\rangle,
  $$
  and that $t\colon X\to\R$ is a homomorphism. We thus compute for all
  $\phi \in \E_0(X)$,
  \begin{eqnarray*}
    \langle\A T * \A S, \phi\rangle &=& \langle\A T , r\mapsto \langle \A S, 
    s\mapsto \phi(r+s)\rangle\rangle \\
    &=& \langle T , x\mapsto \langle S, y\mapsto \phi(t(x)+t(y))
    e^{\rho (t(x)+ t(y))}\rangle\rangle \\
    &=& \langle T , x\mapsto \langle S, y\mapsto \phi(t(xy))e^{\rho (t(xy))}
    \rangle\rangle \\
    &=& \langle\A(T*S), \phi\rangle 
  \end{eqnarray*}
  which proves \eqref{abelalg}. Next, we prove the second claim. Using
  $\langle g,h \rangle = (g*h)(e)$ and commutativity of the
  convolution of radial functions, we have for all $g, \phi \in \D_0(X)$
  \begin{eqnarray*}
    \langle g * af, \phi \rangle &=& \langle af, g * \phi \rangle \\
    &=& \langle f, \A g * \A \phi \rangle \quad \text{(see \eqref{abelalg})}\\
    &=& \langle A g * f, \A \phi \rangle \\
    &=& \langle a (\A g * f),\phi \rangle.
  \end{eqnarray*}
  Now, from the continuity of $a, \A$ and the density of $\D_0(X)$ in
  $\E_0'(X)$ we get \eqref{abelalg2}.\qed
\end{proof}

We will prove in the next paragraph that $a$ is bijective. Assuming this for the moment we have
\begin{prop}\label{inder} The following diagram commutes:
\begin{equation}\begin{CD}
       \E_0'(X) \times \E_0(X) @>*_X>> \E_0(X)\\
       @VV{\A \times \B}V @VV{\B}V \\
       \E_0'(\R) \times \E_0(\R) @>*_\R>> \E_0(\R)
   \end{CD}
\end{equation}
where $$\B:=a^{-1}\colon\E_0(X)\to\E_0(\R) \ . $$
\end{prop}

\subsection{Bijectivity of the dual Abel transform}

In this paragraph we show

\begin{prop}\label{abelinver}
  The maps
  $$a\colon\E_0(\R)\to\E_0(X)\qmbox{and}\A=a'$$ 
  are topological isomorphisms.
\end{prop}

\begin{proof} The proposition is essentially a consequence of the
  bijectivity of $\A\colon\D_0(X)\to\D_0(\R)$ (see \cite{ADY-96}). For
  injectivity of $a$, assume $aw=0$ for some $w\in\E_0(\R)$. Then
  $$\langle u, aw\rangle=\langle\A u, w\rangle=0 $$
  for all $u\in\D_0(X)$. Since $\A$ is surjective it follows that
  $w=0$.  Surjectivity of $a$ follows from the explicit calculation of
  $a^{-1}: \E_0(X) \to \E_0(\R)$. Introducing the bijective map $\Phi:
  C^\infty([1,\infty)) \to \E_0(\R)$,
  \[ (\Phi f)(r) = f(\cosh r), \] 
  one first observes that
  \[
  \Phi^{-1} \frac{d}{d(\cosh r)} \Phi = \frac{d}{dt}, \quad \Phi^{-1}
  \frac{d}{d(\cosh(r/2)} \Phi = 2 \sqrt{2} (t+1)^{1/2} \frac{d}{dt}.
  \]
  Using the explicit formulas for $\A^{-1}: \D_0(\R) \to \D_0(X)$ in
  \cite{ADY-96}, lengthy (but straightforward) computations yields in
  the case $p = 2k$ and $q=2l$:
  \begin{multline*}
    (\Phi^{-1} a^{-1} \Phi u)(t) = C_{p,q} (t+1)^{1/2} (t-1)^{1/2} \times\\
    \left( \frac{d}{dt} (t+1)^{1/2} \right)^k \left( \frac{d}{dt}
    \right)^l (t+1)^{l-1/2} (t-1)^{l+k-1/2} u(t),
  \end{multline*}
  and in the case $p =2k$ and $q=2l-1$:
  \begin{multline*}
    (\Phi^{-1} a^{-1} \Phi u)(t) = C_{p,q} (t+1)^{1/2} (t-1)^{1/2} \times\\
    \left( \frac{d}{dt} (t+1)^{1/2} \right)^k \left( \frac{d}{dt}
    \right)^l (t+1)^{l-1/2} (t-1)^{l+k-1/2}(R_{1/2}^{(k+l-1,l-1)}
    u)(t),
  \end{multline*}
  with suitable constants $C_{p,q}$, and where
  $R_{1/2}^{(\alpha,\beta)}$ is defined in \cite{Ko-84}. Hence $a:
  \E_0(\R) \to \E_0(X)$ is a bijective linear continuous map. By the
  open mapping theorem (see, e.g., \cite[Theorem 17.1]{Tr-67}), $a$ is
  a topological isomorphism. From the Corollary of Proposition
  19.5. in \cite{Tr-67}, we conclude that $\A: \E_0'(X) \to \E_0'(\R)$
  is also a topological isomorphism. \qed
\end{proof}

\subsection{The Paley-Wiener Theorem for the spherical Fourier
  Transform on Distributions}

Let $\bE_0'$ denote the Fourier tansform of the space $\E_0'(\R)$. The
classical Paley-Wiener theorem for distributions (see, e.g.,
\cite{Don-69}) states that $\bE_0'$ consists of all even
entire functions $f: \C \to \C$ of exponential type which are
polynomially bounded on $\R$, i.e., there are constants $C, R \ge 0$
and $m \ge 0$ such that
\[ | f(\lambda) | \le C (1+|\lambda|)^m e^{| {\rm Im}(\lambda) |
  R}. \] We topologize $\bE_0'$ by choosing the subsets $U_a \subset
\bE_0'$ as a fundamental system of neighbourhoods of $0$, where
\[ U_a := \{ f \in \bE_0' \mid | f(\lambda) | \le a(\lambda) \} \] and
$a$ is any continuous positive function of the form
$a(\lambda)=a_1({\rm Re} \lambda)a_2({\rm Im} \lambda)$, where $a_1$
dominates all polynomials and $a_2$ dominates all linear
exponentials. Then the Fourier transform is a topological isomorphism
$\E_0'(\R) \to \bE_0'$, by \cite[Theorem 5.19]{Ehr-70}.

As a direct consequence of this fact, Proposition \ref{abelinver},
formulas \eqref{fourabel} and \eqref{abelalg}, we obtain (see also
\cite[Theorem 4]{FlJ-72}):

\begin{theorem} \label{PaWie} 
  The spherical Fourier transform
  $$\F T(\lambda)=\langle T, \phi_\lambda\rangle $$
  defines a topological isomorphism
  $$\F\colon\E_0'(X)\to\bE_0'.$$
  Furthermore, for distributions $T,S \in \E_0'(X)$, we have
  $$\F (T*S)=\F T\cdot\F S.$$
\end{theorem}


\section{Spherical spectral synthesis in Damek Ricci Spaces}\label{sss}

In this section we prove spherical spectral synthesis in $\E_0(X)$.
We begin with two applications of Proposition \ref{inder}.

\begin{lemma}\label{inderappl1}
  Let $T \in \E_0'(X)$. Then
  \[ \langle T, \varphi_\lambda \rangle = 0 \, \Leftrightarrow\, T*
  \varphi_{\lambda} = 0. \]
\end{lemma}

\begin{proof}
  The implication $T * \varphi_\lambda = 0 \Rightarrow \langle T,
  \varphi_\lambda \rangle = 0$ is obvious. Now, assume $\langle T,
  \varphi_\lambda \rangle = 0$. Using \eqref{abelfilambda} and Proposition
  \ref{inder} we obtain
  \[ \langle\A T,\psi_\lambda\rangle = \langle\A T,
  \B\varphi_\lambda\rangle = \langle T, \varphi_\lambda\rangle = 0. \]
  Using, again, Proposition \ref{inder}, we also obtain
  \[ \B(T * \varphi_\lambda) = (\A T) * (\B\varphi_\lambda) = (\A T) *
  \psi_\lambda. \] 
  Moreover,
  \begin{multline*}
    (\A T * \psi_\lambda)(t) = \langle\A T, s \mapsto
    \psi_\lambda(t-s) \rangle=\\
    \frac{1}{2}(e^{i \lambda t}\langle \A T, \phi_\lambda \rangle +
    e^{-i \lambda t}\langle \A T, \phi_\lambda
    \rangle)=\psi_\lambda(t)\langle\A T, \psi_\lambda \rangle=0.
  \end{multline*}
  Since $\B$ is an isomorphism, we conclude that $T * \varphi_\lambda
  = 0$. \qed
\end{proof}

Recall that a variety $V \subset \E_0(X)$ is a proper closed subspace
satisfying $\E_0'(X) * V \subset V$.

\begin{lemma} \label{inderappl2} Let $V \subset \E_0(X)$ be a
  variety. If $\varphi_{\lambda,k} \in V \backslash \{ 0\}$ then also
  $\varphi_{\lambda,l} \in V$ for all $0 \le l \le k$.
\end{lemma}

\begin{proof}
  Let $W := \B(V)$. By Proposition \ref{inder}, $W \subset \E_0(\R)$
  is also a variety. From \eqref{abelfilambda} we have
  $\psi_{\lambda,k} = \B(\varphi_{\lambda,k})$. So it remains to prove
  that
  \[
  \psi_{\lambda,k} \in W\backslash \{0\} \, \Rightarrow \, \psi_{\lambda,l} \in W
  \quad \forall\, 0 \le l \le k.
  \]
  In the case $\lambda \neq 0$, we restrict our considerations to
  $k=1$ (the case $k \ge 2$ is proved similarly.)  Note that $f \in W$
  implies $f_s + f_{-s} \in W$, where $f_s(t) = f(s+t)$.  Therefore,
  we have
  \[
  (\psi_{\lambda,1})_s + (\psi_{\lambda,1})_{-s} - 2 \cos(s\lambda)
  \psi_{\lambda,1} = -2s \sin(\lambda s) \psi_\lambda(x).
  \]
  Consequently, we have $\psi_\lambda \in W$. If $\lambda =0$,
  $\psi_{0,k} =0$ if $k$ is odd and $\psi_{0,k}$ are monomials if $k$
  is even. Here,
 \[
  (\psi_{0,2})_s + (\psi_{0,2})_{-s} - 2 \psi_{0,2} = -2s^2 \psi_0.
  \]
  (The case $k \ge 4$ is treated similarly.)
  \qed
\end{proof}

Next, we prove an equivalent formulation of {\em spherical spectral synthesis}
in Damek-Ricci spaces (see Theorem \ref{blackydrintro}):

\begin{theorem} \label{blackydr} 
  Let $f \in \E_0(X)$ be a non-zero mean periodic function, i.e.,
  \[ \{ 0 \} \neq V_0^f := \overline{\{ T * f \mid T \in \E_0'(X) \}}
  \neq \E_0(X). \] Then $f$ is the limit of finite linear combinations
  of functions in
  \[ {\rm spec}_0\, f := \{ \varphi_{\lambda,k}\in
  V_0^f \backslash \{0 \} \mid k \in \N_0, \lambda \in \C \}.
  \]
\end{theorem}

\begin{proof}
  Proposition \ref{inder} implies that
  \[ \B(V_0^f) = \B(\overline{\{ T*f \mid T \in \E_0'(X)\}}) =
  \overline{ \{ S*\B f \mid S \in \E_0'(\R)\} } = V_0^{\B(f)}, \] and,
  by \eqref{abelfilambda},
  \[ \B({\rm spec}_0\, V_0^f) = {\rm spec}_0\, V_0^{\B(f)}. \] Now,
  the theorem follows immediately from Theorem \ref{blackyeven}.\qed
\end{proof}

\begin{remark}
  Theorem \ref{blackydrintro} is a direct consequence of the above Theorem
  since every function in a variety $V \subset \E_0(X)$ is mean periodic. 
\end{remark}

\begin{cor}[Spherical Spectral Analysis] \label{key1}
  If $f \in \E_0(X)$ is mean periodic and ${\rm spec}_0\, f$ is empty,
  then $f = 0$.
\end{cor}

\medskip

The following corollary will be used in the next section.

\begin{cor} \label{key} 
  Let $\Pp$ be a non-empty set of distributions in $\E'_0(X)$. Then
  the following two statements are equivalent:
  \begin{itemize}
  \item[a)] There exists a non-zero function $f \in \E_0(X)$ such that
    $T * f = 0$ for all $T \in \Pp$.
  \item[b)] There exists $\lambda \in \C$ such that
    \[ \F T(\lambda) = 0 \quad \forall \, T \in \Pp. \]
\end{itemize} 
\end{cor}

\begin{proof}
  We first prove b) implies a): If there exists $\lambda \in \C$ with
  $\F T(\lambda) = \langle T,\varphi_\lambda\rangle = 0$ for all $T \in \Pp$,
  Lemma \ref{inderappl1} yields that $T * \varphi_\lambda = 0$ for all $T
  \in \Pp$. Thus, a) is satisfied with $f=\varphi_\lambda$.

  \medskip

  In the proof of a) implies b) we assume that $\Pp \neq \{ 0 \}$.
  (The case $\Pp = \{ 0 \}$ is trivial.) Since $\Pp$ contains at least
  one non trivial distribution, we have $V_0^f \neq \E_0(X)$.  From
  Corollary \ref{key1} we conclude that ${\rm spec}_0 f \neq
  \emptyset$. Then there exists a non-zero $\varphi_{\lambda,k} \in
  V_0^f$. Using Lemma \ref{inderappl2} we conclude that
  $\varphi_\lambda \in V_0^f$, which implies that $T * \varphi_\lambda
  = 0$ for all $T \in \Pp$.  Consequently, we have
  \[ \F T(\lambda) = \langle T, \varphi_\lambda \rangle = (T * \check
  \varphi_\lambda)(e) = 0 \] 
  for all $T \in \Pp$. \qed
\end{proof}


\section{Applications: Two-radius theorems in Damek Ricci Spaces}
\label{tworad}

The following lemma is needed for the proofs of the two-radius theorems.

\begin{lemma} \label{TcheckTpif} Let $T \in \E_0'(X)$, $f \in
  \E(X)$. If $T * \check f = 0$ then $T * (\pi f) = 0$.
\end{lemma}

\begin{proof}
  Let $\rho_\epsilon \in \D_0(X)$ be a Dirac sequence. Then $T *
  \rho_\epsilon =: g_\epsilon \to T$, as $\epsilon \to 0$. Using
  \eqref{piprop} and $f * g = (\check g * \check f)^\vee$, we conclude that
  \begin{multline*}
    T * (\pi f) = \lim_{\epsilon \to 0} g_\epsilon * (\pi f) =
    \lim_{\epsilon \to 0}
    (\pi f) * g_\epsilon = \lim_{\epsilon \to 0} \pi( f * g_\epsilon) =\\
    \lim_{\epsilon \to 0} \pi((g_\epsilon * \check f)^\vee) = \pi( (T*
    \check f)^\vee ) = 0.
  \end{multline*}
  \qed
\end{proof}

\noindent
\PPP
It suffices to prove the theorems for smooth functions only. This
is because the averaging operators are continuous
with respect to uniform convergence on compacta and $\E(X)\subset
C(X)$ is dense.

The proof proceeds by contradiction:
Let $r_1, r_2 > 0$ avoid the set described in the theorem and
assume that
\begin{equation} \label{VV} 
  \V = \{ f \in \E(X) \mid \langle T_{r_1},
  f_x \rangle = 0, \langle T_{r_2}, f_x \rangle = 0 \, \forall \, x
  \in X \} \neq \{ 0 \},
\end{equation}
where the corresponding families of distributions are in each case
\begin{eqnarray} 
\label{Trth1}\langle T_r, f \rangle &=& \int_{B_r} f, \\
\label{Trth2}\langle T_r, f \rangle &=& \int_{S_r} f. 
\end{eqnarray}
Obviously, $\V$ is invariant under left-translations (isometries) in
$X$ (i.e., $f \in \V \Rightarrow f_x \in \V$ for all $x \in
X$). Therefore we can find a function $f \in \V$ with $f(e) \neq
0$. Then $\pi f \neq 0$ and, using $\langle T, f_x \rangle = (T *
\check f)(x^{-1})$, Lemma \ref{TcheckTpif} shows that $T_{r_1}*\pi f =
T_{r_2} * \pi f = 0$.

Now, for $\Pp = \{ T_{r_1}, T_{r_2}\}$, Corollary \ref{key} implies
that there exists a $\lambda \in \C$ with
$$ \F T_{r_1}(\lambda) = \F T_{r_2}(\lambda) = 0. $$
By the following Lemma, we obtain a contradiction to the choice of
the radii $r_1,r_2$ at the beginning of the proof. \qed

\begin{lemma} \label{explicitfour} 
  Let $X = X^{(p,q)}$ be a Damek-Ricci space of dimension $n+1 =
  p+q+1$.
  \begin{itemize}
  \item[a)] Let $T_r$ be defined as in \eqref{Trth1}. Then
    \[ \F T_r(\lambda) = \frac{2^n \pi^{n/2}}{\Gamma(1+n/2)} \left(
      \sinh(r/2) \right)^n \left( \cosh(r/2) \right)^{q-1}
    \varphi_\lambda^{(p,q+2)}(r). \]
  \item[b)] Let $T_r$ be defined as in \eqref{Trth2}. Then
    \begin{equation}\label{ftsppro} \F T_r(\lambda) = \langle T_r,
      \varphi_\lambda \rangle
      = \vol(S_r) \, \varphi_\lambda^{(p,q)}(r). \end{equation}
  \end{itemize}
\end{lemma}

\begin{proof}
  b) is obvious. For the proof of a), note that
  \[ \varphi_\lambda^{(p,q)}(r) = F(\rho-i \lambda, \rho + i \lambda,
  n/2, - \sinh^2(r/2) ) \] 
  with $\rho= p/4 + q/2$. Choosing
  \[ z = -\sinh^2{r/2},\ a = \frac{p}{4}+ \frac{q}{2} - i \lambda,\ b
  = \frac{p}{4}+ \frac{q}{2} + i \lambda,\ c = \frac{n}{2}, \] 
  and using (see \cite[Formula 15.2.9]{AS-72})
  \begin{multline*}
    \frac{d}{dz}\left(z^c (1-z)^{a+b+1-c} F(a+1,b+1,c+1,z) \right) =\\
    c z^{c-1}(1-z)^{a+b-c} F(a,b,c,z),
  \end{multline*}
  a straightforward calculation yields result a). \qed
\end{proof}

\medskip

\noindent
\PPPP 
Let $r_1, r_2 > 0$ avoid the set described in the theorem,
$T_r$ be defined by
$$ \langle T_r, f \rangle = \frac{1}{\vol(S_r)}\left( \int_{S_r}f \right)
- f(e), $$ 
and $\V$ be as in \eqref{VV}. Observe that
$\ker\Delta\subset\V$ and that $\ker\Delta\cap\E_0(X)$ is spanned by
the constant function $\varphi_{i\rho}=\varphi_{- i\rho}= 1$. Assume
there is a function $f \in \V$ with $\Delta f \neq 0$. Since $\V$ and
$\Delta$ are invariant under left-translations, we can assume that
$\Delta f(e) \neq 0$. Let $g=\pi f$. Since $\Delta$ and $\pi$ commute,
we have $\Delta g(e) \neq 0$ and Lemma \ref{TcheckTpif} implies that
$g \in \V \cap \E_0(X)$. Note that $g$ is a non-zero mean periodic
function (since all functions $h$ in $V_0^g$ satisfy $T_{r_j}*h =
0$ and thus $h(r_j)=h(0)$ for $j=1,2$). 

Next, we show that ${\rm spec}_0 g = \{ 1 \}$. Let
$\varphi_{\lambda,k}\in V_0^g \backslash \{ 0 \}$. We will show that $\lambda =
\pm i\rho$ and $k=0$. By Lemma \ref{inderappl2}, we also have
$\varphi_\lambda \in V_0^g$ and, therefore,
$\varphi_\lambda(r_1)=\varphi_\lambda(r_2)=\varphi_\lambda(0)=1$.
This implies that $\lambda = \pm i \rho$. If $k \ge 1$, then we also
must have $\varphi_{i \rho,1} \in V_0^g$ and thus
$\varphi_{i\rho,1}(r_j) = \varphi_{i\rho,1}(0)$. We have
$\varphi_{i\rho,1}=a\psi_{i\rho,1}$ and
$$\psi_{i\rho,1}(r)=\left.\frac{d}{d\lambda}\right|_{\lambda=i\rho}
\cos(\lambda r)=- r\sin(i\rho r)=- ir\sinh(\rho r)$$ 
is $(-i)$ times a positive function for $r>0$. Since the dual Abel
transform is multiplication by a positive real function followed by an
averaging operator, it preserves positivity. Thus
$\varphi_{i\rho,1}(r)$ doesn't vanish for $r > 0$. But
$\varphi_{i\rho,1}(e) = (a \psi_{i\rho,1})(0) =0$ and, consequently,
we must have $k=0$.

By Theorem \ref{blackydr}, $g$ is a constant function, contradicting to
$\Delta g \neq 0$. \qed

\end{document}